\documentclass{elsarticle}

\usepackage{amsmath,amssymb,amsthm}
\usepackage{mathptmx}
\usepackage{epsfig}
\usepackage{verbatim}
\usepackage{float}
\usepackage{subfig}
\usepackage{multirow,bm,url}
\usepackage{algorithm,algpseudocode}
\graphicspath{ {./Figures/} }

\begin{document}

\begin{frontmatter}

\title{Adaptive Discontinuous Galerkin Finite Elements for Advective Allen-Cahn Equation}

\author[mu]{Murat Uzunca\corref{cor}}
\author[ay]{Ay\c{s}e Sar{\i}ayd{\i}n-Filibelio\u{g}lu}

\address[mu]{Department of Mathematics, Sinop University, 57000, Sinop, Turkey}
\address[ay]{Turkish Scientific and Technological Research Council (T\"UB{\.I}TAK), 06100, Ankara, Turkey}
\cortext[cor]{Corresponding author.}

\journal{Numerical Algebra, Control \& Optimization (Accepted)}

 \begin{abstract}
We apply a space adaptive interior penalty discontinuous Galerkin method for solving  advective Allen-Cahn equation with expanding and contracting velocity fields. The advective Allen-Cahn equation is first discretized in time and the resulting semi-linear elliptic PDE is solved by an adaptive algorithm using a residual-based a posteriori error estimator. The a posteriori error estimator contains additional terms due to the non-divergence-free velocity field. Numerical examples demonstrate the effectiveness and accuracy of the adaptive approach by resolving the sharp layers accurately.
\bigskip

\noindent MSC 2010: 65M20; 65M22; 65M50; 65M60
 \end{abstract}
 \begin{keyword}
Advective Allen-Cahn equation; discontinuous Galerkin method; Rothe's method; adaptivity.
\end{keyword}

 \end{frontmatter}

\section{Introduction}

Interfacial dynamics has great importance in the modeling of multi-phase flow and it plays an important role in different scientific and industrial applications such as micro-structure evolution and grain growth in material science \cite{chen02pfm}, binary fluids flow movement \cite{osher88fpc}, and complex interfacial dynamics \cite{hohenberg77tdc}.

There have been various diffuse interface models for multi-phase flow \cite{khan09pfm,liu12dis}. In this study, we consider a specific model of diffuse interface for two phase flow; Allen-Cahn (AC) equation with advection. AC equation without advection is the most known dynamical model for diffuse interface dynamics. In the past it was investigated extensively with many numerical and analytical methods.    But few studies are dealing with the numerical solution of the advective AC equation \cite{Shen16,Shen10,Vasconcelos14,Zhang07}. 
The AC model with advection is given by
\begin{subequations}\label{advac}
\begin{align}
 \frac{\partial u}{\partial t} + \nabla \cdot (\bm{V} u) &= \epsilon\Delta u - \frac{1}{\epsilon}f(u), & \hbox{in }  \Omega \times (0,T], \\
u(\bm{x},0) &= g(\bm{x}), & \hbox{in }  \Omega \times \{0\} , \\
	\frac{\partial u}{\partial n} &=0, & \hbox{on } \partial \Omega \times (0,T].
\end{align}
\end{subequations}
In \eqref{advac}, the bi-stable cubic nonlinearity $f(u)$ is such that $f(u)=F'(u)=2u(1 - u)(1 -2u)$ for a double-well potential $F(u)$, and $\bm{V}=(V_1, V_2)^T$ is the prescribed velocity field. The velocity field is either given or computed from the Navier-Stokes equation  \cite{khan09pfm,liu12dis}.
In many cases the velocity field $\bm{V}$ is divergence-free  \cite{Shen16,Shen10} and satisfies incompressible Navier-Stokes equations. In this study we consider non-divergence-free velocity fields $\nabla \cdot\bm{V} \neq 0$, either expanding  $\nabla \cdot \bm{V}> 0$ or contracting  $\nabla \cdot \bm{V} < 0$ \cite{liu12dis}. In other words we consider advective AC equation in compressible fluids.

Problems with  surface tension in two-phase fluids are known as multi-scale problems with  two different time scales, the small surface tension, and the convection time scale, which results in computational stiffness. Actually, there exists  three main algorithms: the sharp interface algorithm method, the level-set algorithm method and the diffuse interface method \cite{liu12dis}. The numerical simulations are illustrated using finite elements method in space and semi-implicit schemes or semi-implicit schemes with splitting in time \cite{barrett98fea,liu11dsm}.

Numerical solutions of advective AC equation may exhibit unphysical oscillations at the interior layers due to convection and sharp fronts may occur due to the non-linear reaction term. Since the standard FEMs are known to produce strong oscillations around layers, adaptive algorithms are developed to tackle the unphysical oscillations and shocks. By refining the mesh locally at layers and sharp fronts, accurate solutions are obtained with less degrees of freedom (DoFs) and computational time. The major part of an adaptive algorithm is the estimation of the local errors and refine the elements with large estimated errors. A posteriori error estimation is the main tool to estimate the local errors which uses the approximate solution and the given problem data \cite{ainsworth00pee,babuska01fem,verfurth05rpes}. Since the discontinuous Galerkin (dG) methods have the flexibility on adaptive meshes, a  posteriori error estimators are developed using dG discretization \cite{cangiani13adg,ern10grd,hoppe08caa,houston02dhp,karakashian07cad,riviere03pee}.

We develop an adaptive strategy for the numerical solution of advective AC equation \eqref{advac} for contracting and expanding flows. Because the solutions do not show strong variations with respect to time, only space adaptivity is applied.  Usually the time dependent PDEs are discretized first in space and integrated in time. Here we discretized first in time by Rothe's method \cite{deuflhard12ans}. Then, the resulting semi-linear elliptic equations are solved with the adaptive version of symmetric interior penalty Galerkin (SIPG) method using upwinding for the convective term. The adaptive strategy is based on a residual based a posteriori error estimation \cite{cangiani13adg,uzunca14adg}. In \cite{uzunca14adg} we developed a posteriori error bounds with respect to the energy norm induced by the SIPG formulation of the system given  for semi-linear diffusion-convection-reaction equations with divergence-free velocity field. In this study we develop a posteriori error estimates for the advective AC equation with non-divergence-free convective terms. The extra terms in the error estimates coming from the non-divergence-free vector field are added to the reaction terms.  Numerical tests show that the proposed adaptive algorithm can efficiently detect  the layers and accurate solutions around these layers can be obtained.

The paper is organized as follows. We first derive the space/time discretization of the advective AC equation \eqref{advac} in Section~\ref{sec:rothe}. A detailed explanation for space adaptivity algorithm is presented in Section~\ref{sec:adap}, where we derive a posteriori error bounds for a stationary problem. Finally, we present two numerical examples with expanding and contracting flows in Section~\ref{sec:numerics} to demonstrate the performance of the adaptive approach.

\section{Discretization (Rothe's method)}
\label{sec:rothe}

Since the solutions of advective AC equation \eqref{advac} do not change much with the evolution of time, we apply only space adaptivity ($h$-adaptivity),  which utilizes a posteriori error estimation or stationary (elliptic) problems at each time step. We first discretize \eqref{advac} in time for obtaining by Rothe's method \cite{deuflhard12ans}. Then, on each time interval, space discretization and adaptivity is applied to the elliptic stationary problem.

\subsection{Time discretization}

For the semi-discrete scheme, we consider the uniform partition $0=t_0<t_1<\ldots < t_J=T$ of the time interval $[0,T]$ with the uniform time step-size $\tau=T/J$, and with the time intervals $I_{k} = (t_{k-1},t_{k}]$. For $k=1,\ldots , J$, we denote by $u^k(\bm{x})$ the approximate solution at the time instance $t=t_k$, i.e. $u^k(\bm{x})\approx u(\bm{x},t_k)$, and we set $u^0(\bm{x}) = g(\bm{x})$. We use backward Euler as the time integrator to discretize the advective AC \eqref{advac}, for $k=1,\ldots , J$, given $u^{k-1}(\bm{x})$ find $u^{k}(\bm{x})$ satisfying
\begin{equation}\label{ellipticAC0}
\begin{split}
\frac{u^{k}(\bm{x}) - u^{k-1}(\bm{x})}{\tau} &- \epsilon\Delta u^{k}(\bm{x})  + \bm{V}\cdot \nabla u^{k}(\bm{x}) \\
 &+ \; (\nabla \cdot \bm{V}) u^{k}(\bm{x}) + \frac{1}{\epsilon}f(u^{k}(\bm{x}) ) = 0.
\end{split}
\end{equation}
For each $k=1,2,\ldots , J$, the equation \eqref{ellipticAC0} can be written in the form of a semi-linear elliptic problem as
\begin{equation}\label{ellipticAC}
\alpha u^{k} - \epsilon \Delta u^{k} + \bm{V}\cdot \nabla u^{k} + r(u^{k})=l(u^{k-1}),
\end{equation}
where
\begin{subequations}\label{alpha}
\begin{align}
\alpha := \alpha (\bm{x}) &= \frac{1}{\tau} + \nabla \cdot \bm{V}(\bm{x}) , \\
r(u^{k}) &= \frac{1}{\epsilon}f(u^{k}(\bm{x}) ), \\
l(u^{k-1}) &= \frac{1}{\tau} u^{k-1}.
\end{align}
\end{subequations}

We assume that the non-linear reaction term is bounded and locally Lipschitz continuous, i.e., satisfy for any $s, s_1, s_2\ge 0$, $s,s_1, s_2 \in \mathbb{R}$ the following conditions
\begin{subequations}\label{nonlas}
\begin{align}
|r(s)| &\leq C , \quad C>0   \\
\| r(s_1)-r(s_2)\|_{L^2(\Omega)} &\leq  L\| s_1-s_2\|_{L^2(\Omega)}, \quad L>0.
\end{align}
\end{subequations}

Moreover, we assume that there is a non-negative constant $\kappa_0$ satisfying
\begin{subequations}\label{erassmp}
\begin{align}
\alpha (\bm{x})-\frac{1}{2}\nabla\cdot \bm{V} (\bm{x}) &\geq \kappa_0, \label{erassmp1} \\
\| -\nabla\cdot \bm{V}(\bm{x}) +\alpha (\bm{x})\|_{L^{\infty}(\Omega)} &\leq c^*\kappa_0, \label{erassmp2}
\end{align}
\end{subequations}
for a positive constant $c^*$. The identity \eqref{erassmp1} guarantees the coercivity of the bilinear form $a_h(\cdot , \cdot)$ in \eqref{dgAC}, and \eqref{erassmp2} is needed to prove the reliability of a posteriori error estimator in Section~\ref{sec:adap}. According to our setting for $\alpha$ in \eqref{alpha}, the assumptions in \eqref{erassmp} are equivalent that
\begin{subequations}\label{erassmpT}
\begin{align}
\frac{1}{\tau}+\frac{1}{2}\nabla\cdot \bm{V} (\bm{x}) &\geq \kappa_0, \\
\| \tau^{-1}\|_{L^{\infty}(\Omega)} &\leq c^*\kappa_0.
\end{align}
\end{subequations}

\subsection{Space discretization}

For the space discretization, we use symmetric interior penalty Galerkin (SIPG) method \cite{arnold82ipf,riviere08dgm}, which is a member of the family of discontinuous Galerkin (dG) finite elements methods, and we apply upwinding \cite{lesaint74fes,reed73tmm} for the convective term. The dG methods exhibit attractive properties of both classical finite elements method (FEM) and finite volume method (FVM). The functions in dG spaces are discontinuous along the inter-element boundaries, which makes dG methods flexible. In addition, different  order of basis functions on each element can be used within dG schemes ($p$-adaptivity). Hence, it allows to use {\it hp}-adaptivity methods \cite{solin03hof} which arranges the mesh elements and also the order of polynomials on each element adaptively. Further, the dG methods locally conserve several physical quantities such as mass and energy, which plays an important role in the flow and transport problems. Moreover, the sharp gradients or the singularities in the mesh can be locally detected by the fully discontinuous polynomial representation of the solution in dG schemes. A sample archive of MATLAB codes together with a dGFEM tutorial can be found at \url{https://github.com/muzunca/dGFEM-Tutorial.git}.

On the $k$-th time interval, let $\{\mathcal{T}_h^k\}_h$ be a family of shape regular meshes with triangular elements $\{E_i\}$ such that $\overline{\Omega}=\cup \overline{E}_i$ and $E_i\cap E_j=\emptyset$ for $E_i$, $E_j$ $\in\mathcal{T}_h^k$. By shape regularity, we mean that there exists a constant $c_0$ such that
\[
\max_{E \in \mathcal{T}_h^k}\frac{h^{2}_{E}}{\left|E\right|}\leq c_0,
\]
where $h_E$ is the diameter and ${\left|E\right|}$ is the area of the triangular element $E$.  We split the set of all edges $\mathcal{F}_h^k$ into the set of interior edges $\mathcal{F}^{k}_{h,0}$ and the set of boundary edges $\mathcal{F}^{k}_{h,\partial}$ so that $\mathcal{F}_{h}^k=\mathcal{F}^{k}_{h,0} \cup \mathcal{F}^{k}_{h,\partial}$. We set the finite dimensional solution and test function space by
\[
V_h^k=\left\{ \upsilon\in L^2(\Omega ) : \upsilon|_{E}\in\mathbb{P}^q(E) ,\; \forall E\in \mathcal{T}_h^k \right\},
\]
where $\mathbb{P}^q(E)$ denotes the set of all polynomials of degree at most $q$ on the element $E\in\mathcal{T}_h^k$. Note that the functions in $V_h^k \not\subset H^{1}_{0}(\Omega)$ are allowed to be discontinuous along the inter-element boundaries, thus, in contrast to continuous FEM, the dG methods are suitable to use a non-conforming space.

The discontinuities of the functions in $V_h^k$ along the inter-element boundaries lead to different traces from the neighboring elements sharing an edge $e$. Let the edge $e$ be a common edge for two elements $E_i$ and $E_{j}$ (w.l.o.g. assume $i<j$). Then for a scalar function $\upsilon\in V_h^k$, there are two traces of $\upsilon$ along $e$, denoted by $\upsilon_{|_{E_i}}$ from inside $E_i$ and  $\upsilon_{|_{E_j}}$ from inside $E_j$. Then, the jump and average of $\upsilon$ across the edge $e$ are defined as
\[
[\upsilon]_e= \upsilon_{|_{E_i}}\bm{n}_e- \upsilon_{|_{E_j}}\bm{n}_e , \quad \{ \upsilon\}_e=\frac{1}{2}(\upsilon_{|_{E_i}}+ \upsilon_{|_{E_j}}),
\]
where $\bm{n}_e$ is the unit normal to the edge $e$ oriented exterior to $E_i$. Similarly, we set the jump and average values of a vector field $\nabla \upsilon$ on e as
\[
\begin{split}
[\nabla \upsilon]_e &= \nabla \upsilon_{|_{E_i}}\cdot \bm{n}_e- \nabla \upsilon_{|_{E_j}}\cdot \bm{n}_e , \\
\{ \nabla \upsilon\}_e &=\frac{1}{2}(\nabla \upsilon_{|_{E_i}}+ \nabla \upsilon_{|_{E_j}}).
\end{split}
\]
For a boundary edge $e\subset E_i\cap \partial\Omega$, we set
\[
\begin{split}
[\upsilon]_e= \upsilon_{|_{E_i}}\bm{n} , \; &\{ \upsilon\}_e=\upsilon_{|_{E_i}}, \\
 [\nabla \upsilon]_e=\nabla \upsilon_{|_{E_i}}\cdot \bm{n}, \; &\{ \nabla \upsilon\}_e=\nabla \upsilon_{|_{E_i}},
\end{split}
\]
where $\bm{n}$ is the unit outward normal to the boundary at $e$. We further define the sets of inflow and outflow boundary parts, respectively, by
\[
\Gamma^{-}_{h} = \left\{\bm{x} \in \partial \Omega: \bm{V}(\bm{x})\cdot \bm{n}(\bm{x})<0\right\}, \quad \Gamma^{+}_{h} =\partial\Omega\setminus \Gamma^{-}_{h}.
\]
The set of inflow and outflow boundary edge parts for an element $E$ is defined in a similar way by
\[
\partial E^{-}_{h} = \left\{\bm{x} \in \partial E: \bm{V}(\bm{x})\cdot \bm{n}_e(\bm{x})<0\right\}, \; \partial E^{+}_{h} = \partial E\setminus\partial E^{-}_{h},
\]
where $\bm{n}_e(\bm{x})$ is the unit outward normal vector to the element boundary $\partial E$ at $\bm{x}$. Moreover, for an interior edge $e$, we denote the trace of a scalar function $\upsilon$ from inside the element $E$ by $\upsilon^{in}$ and from outside the element $E$ by $\upsilon^{out}$.

We multiply the continuous (in space) equation \eqref{ellipticAC} by a test function $\upsilon_h^k\in V_h^k$, integrate over $\Omega$, apply Green's identity for diffusive term together with the upwinding for the convective term, and we obtain the following weak problem: set $u_h^{0}$ be the $L^2$-projection of $u^{0}(\bm{x})$ onto the dG space $V_h^k$, for $k=1,2,\ldots , J$, given $u_h^{k-1}$, $\forall \upsilon_{h}^k \in V_h^k$ find $u_h^{k}\in V_h^k$ satisfying
\begin{equation}\label{dgAC}
a_{h}(u_h^{k},\upsilon_{h}^k) + r_h(u_h^{k}) = l_{h}^{k-1}(\upsilon_{h}^k),
\end{equation}
where the dG bilinear form is such that
\[
a_{h}(u,\upsilon):= D_{h}(u,\upsilon) + O_{h}(u,\upsilon) + K_{h}(u,\upsilon) + J_{h}(u,\upsilon),
\]
where the bilinear terms are given by
\begin{equation}\label{dohj}
\begin{split}
D_{h}(u,\upsilon)= & \sum_{E \in \mathcal{T}_h^k}\int_E \left( \epsilon \nabla u \cdot\nabla \upsilon + (\alpha - \nabla\cdot\bm{V})u\upsilon \right)d\bm{x}, \\
O_{h}(u,\upsilon)= & - \sum \limits_{E \in {\mathcal{T}}_{h}^k} \int_{E} \bm{V}u\cdot \nabla \upsilon d\bm{x} + \sum \limits_{E \in {\mathcal{T}}_{h}^k} \int_{\partial E_h^+\cap \Gamma_h^{+}} \bm{V}\cdot \bm{n}_E u \upsilon ds \\
& + \sum \limits_{E \in {\mathcal{T}}_{h}^k}\int_{\partial E_h^+\setminus\partial\Omega } \bm{V}\cdot \bm{n}_E u(\upsilon -\upsilon^{out}) ds, \\
K_{h}(u,\upsilon)= & - \sum_{e\in \mathcal{F}^{k}_{h,0}}\int_e\left( \{\epsilon\nabla u\}[\upsilon]  +  \{\epsilon\nabla \upsilon\}[u] \right)ds, \\
J_{h}(u,\upsilon)= & \sum_{e\in \mathcal{F}^{k}_{h,0}} \frac{\sigma \epsilon}{h_{e}}\int_e [u][\upsilon]ds,
\end{split}
\end{equation}
with $\sigma$ denoting the penalty parameter which should be sufficiently large for the stability of the dG scheme \cite{riviere08dgm}.
The nonlinear form $r_h(u)$ and the linear right hand side $l_h^{k-1}(\upsilon)$ in \eqref{dgAC}  are given by
\[
\begin{split}
r_h(u) &= \sum_{E \in \mathcal{T}_h^k}\int_E r(u)\upsilon d\bm{x}, \\
l_h^{k-1}(\upsilon) &= \sum_{E \in \mathcal{T}_h^k}\int_E \frac{1}{\tau} u_h^{k-1} \upsilon d\bm{x}.
\end{split}
\]

\section{Space adaptivity}
\label{sec:adap}

In this section, a space-adaptive procedure is constructed for the SIPG discretized system \eqref{dgAC} of the stationary problem \eqref{ellipticAC} which mimics the $k$-th time step of the semi-discretized advective AC equation \eqref{advac}. We use a residual-based a posteriori error approach to apply the adaptivity in which not only we employ the refinement but also we consider the coarsening phenomena. In \cite{schotzau09rae}, an adaptive scheme using a posteriori error estimates is constructed for stationary linear diffusion-convection-reaction equation using dG. It is extended to stationary diffusion-convection equations with non-linear reaction using dG in \cite{uzunca14adg}, and to parabolic diffusion-convection equations with non-linear reaction in \cite{cangiani13adg,uzunca14tsa}.

The adaptive algorithm starts with a sufficiently coarse initial mesh $\mathcal{T}_h^0$ together with an initial vector $u_h^0$ which is the $L^2$-projection of the initial condition $u^0(\bm{x}):=g(\bm{x})$ onto the initial solution space $V_h^0$.
Then, on each time interval $I_k=(t_{k-1},t_k]$, we consider a single stationary problem \eqref{ellipticAC} with its SIPG discretized system \eqref{dgAC}. Firstly, we solve the discrete system \eqref{dgAC} for the solution $u_h^k$ on the space $V_h^{k-1}$ given that $u_h^{k-1}\in V_h^{k-1}$ is the known solution vector from the previous time interval. Then, it follows the estimation step providing information about the elements for refinement/coarsening. In the estimation procedure, we utilize a posteriori error estimation. An adaptive scheme with the use of an a posteriori error estimation requires two crucial ingredients. One is an error indicator to compute the local errors on each element, and the other is a compatible norm to measure the error. Here, in order to measure the error, we use the following so-called dG norm
\begin{equation}\label{energy_norm}
\|\upsilon\|_{dG}^2 :=  ||| \upsilon |||^2  +  |\upsilon|_C^2,
\end{equation}
which is composed of the energy-like norm
\[
\begin{aligned}
||| \upsilon |||^2 =& \sum \limits_{E \in \mathcal{T}_h^k}\left( \| \epsilon\nabla \upsilon\|_{L^2(E)}^2 + \kappa_0 \| \upsilon\|_{L^2(E)}^2 \right) \\
& + \sum_{e\in \mathcal{F}^{k}_{h,0}} \frac{\sigma \epsilon}{h_{e}}\| [\upsilon]\|_{L^2(e)}^2,
\end{aligned}
\]
and the semi-norm
\[
|\upsilon|_C^2=|\bm{V}\upsilon|_*^2+\sum_{e\in \mathcal{F}^{k}_{h,0}}(\kappa_0 h_e+ \frac{h_e}{\epsilon})\| [\upsilon]\|_{L^2(e)}^2,
\]
where
\[
|\bm{V}\upsilon|_*=\mathop{\text{sup}}_{w\in H_0^1(\Omega )\setminus \{ 0\}}\frac{\int_{\Omega}\bm{V}\upsilon\cdot \nabla w dx}{||| w|||},
\]
and the parameter $\kappa_0$ is the lower bound in \eqref{erassmpT}.
As the indicator on the mesh $\mathcal{T}_h^{k-1}$, we use the following a posteriori error indicator \cite{uzunca14tsa,schotzau09rae,uzunca14adg}:
\begin{align}\label{res}
(\eta_E^k)^2 &= (\eta_{E,R}^k)^2  + (\eta_{E,0}^k)^2, & \forall E\in\mathcal{T}_h^{k-1},
\end{align}
where $\eta_{E,R}^k$ stands for the volume residual on the element $E$, given by
\begin{equation}\label{res1}
\begin{split}
(\eta_{E,R}^k)^2 = & \; \rho_E^2\| l(u_h^{k-1})-\alpha_h u_h^k + \epsilon\Delta_h u_h^k \\
& - \bm{V}_h\cdot\nabla_h u_h^k - r(u_h^k)\|_{L^2(E)}^2,
\end{split}
\end{equation}
while $\eta_{E,0}^k$ denotes the edge residuals coming from the jump of the numerical solution on the interior edges, given by
\begin{equation}\label{res2}
\begin{split}
(\eta_{E,0}^k)^2 = & \; \frac{1}{2}\sum \limits_{e \in \partial E\cap \mathcal{F}_{h,0}^{k-1}}\left\{\epsilon^{-\frac{1}{2}}\rho_e\| [\epsilon\nabla_h u_h^k]\|_{L^2(e)}^2 \right.\\
 & \left. + \left( \frac{\epsilon\sigma}{h_e} + \kappa_0 h_e+\frac{h_e}{\epsilon}\right)\| [u_h^k]\|_{L^2(e)}^2\right\}.
\end{split}
\end{equation}
In \eqref{res1} and \eqref{res2}, $\Delta_h$ and $\nabla_h$ denote the discrete versions of the Laplace and gradient operators, respectively, and $\alpha_h$ and $\bm{V}_h$ are the $L^2$-projections of the data $\alpha$ and $\bm{V}$ onto $V_h^{k-1}$, respectively. The positive weights $\rho_E$ and $\rho_e$ are defined by
\[
\rho_{E}=\min\{h_{E}\epsilon^{-\frac{1}{2}}, \kappa_0^{-\frac{1}{2}}\}, \; \rho_{e}=\min\{h_{e}\epsilon^{-\frac{1}{2}}, \kappa_0^{-\frac{1}{2}}\},
\]
for $\kappa_0 \neq 0$, or $\rho_{E}=h_{E}\epsilon^{-\frac{1}{2}}$ and $\rho_{e}=h_{e}\epsilon^{-\frac{1}{2}}$ when $\kappa_0 =0$. Then, the global error indicator is set as
\begin{equation} \label{errind}
\eta^k=\left( \sum \limits_{E\in\mathcal{T}_h^{k-1}}(\eta_E^k)^2\right)^{1/2}.
\end{equation}
We also introduce the data approximation terms,
\[
(\Theta_E^k)^2 =\rho_E^2(\| \alpha-\alpha_h\|_{L^2(E)}^2 + \| (\bm{V} -\bm{V}_h)\cdot\nabla u_h^k\|_{L^2(E)}^2),
\]
and  the data approximation error
\begin{equation} \label{data}
\Theta^k=\left( \sum \limits_{E\in\mathcal{T}_h^{k-1}}(\Theta_E^k)^2\right)^{1/2}.
\end{equation}
Then, using the definitions in \eqref{errind} and \eqref{data}, we can obtain the a posteriori error bounds
\begin{align}
\| u^k-u_h^k\|_{dG} &\lesssim \eta^k + \Theta^k,  &(\text{reliability}) \label{rel} \\[0.2cm]
\eta^k &\lesssim \| u^k-u_h^k\|_{dG} + \Theta^k,  &(\text{efficiency}) \label{eff}
\end{align}
where $\lesssim$ mimics a bound up to a constant. The proof of the a posteriori error bounds in \eqref{rel} and \eqref{eff} proceeds as the following: first the dG solution $u_h^k$ is rewritten as the direct sum of its conforming part $u_{h,c}^k\in H_0^1(\Omega)\cap V_h^{k-1}$ and the non-conforming remainder term  $u_{h,r}^k\in V_h^{k-1}$:
\[
u_h^k := u_{h,c}^k + u_{h,r}^k,
\]
and from triangle inequality,
\[
\| u^k-u_h^k\|_{dG} \leq \| u^k-u_{h,c}^k\|_{dG} + \| u_{h,r}^k\|_{dG},
\]
where the term $\| u^k-u_{h,c}^k\|_{dG}$ is now well-defined. Then, using coercivity of the bilinear form $a_h$, continuity of the bilinear forms $D_h$, $O_h$, $K_h$ and $J_h$, inf-sup condition, error bounds for approximation and interpolation operator, the assumptions \eqref{nonlas} on the nonlinearity, we derive the following bounds
\[
\begin{split}
\| u_{h,r}^k\|_{dG} &\lesssim \eta^k  ,\\
\| u^k-u_{h,c}^k\|_{dG} &\lesssim \eta^k +  \Theta^k,
\end{split}
\]
which finishes the proof of reliability. For further details on the proof of reliability and efficiency of the a posteriori error estimation we refer to \cite{cangiani13adg,uzunca14tsa,schotzau09rae,uzunca14adg}.

After collecting information about the elements in the estimation step, we mark the elements for refinement/coarsening. We prescribe two tolerances $\bf{stol^r}$ and $\bf{stol^c}$ related to the refinement and coarsening, respectively. We refine an element $E$ for which the corresponding error indicator $\eta_E^k$ is greater than $\bf{stol^r}$, and we coarsen an element $E$ if it is less than $\bf{stol^c}$, provided that $E$ is not an element of the initial mesh. Accordingly, we form the sets $M_R$ and $M_C$ of the elements in $\mathcal{T}_{h}^{k-1}$ to be refined and coarsened defined, respectively, by
\[
\begin{split}
M_R &= \{ E\in \mathcal{T}_{h}^{k-1}\; : \; (\eta_E^k)^2 > \bf{stol^r} \}, \\
M_C &= \{ E\in \mathcal{T}_{h}^{k-1}\; : \; (\eta_E^k)^2 < \bf{stol^c} \}.
\end{split}
\]
Then, we create the new mesh $\mathcal{T}_{h}^{k}$ by refining the elements $E\in M_R$ using the newest vertex bisection method \cite{chen08fem}, and by coarsening the elements $E\in M_C$. As the final stage, we resolve the discrete system \eqref{dgAC} for $u_h^k\in V_h^k$ on the new mesh $\mathcal{T}_{h}^{k}$. However, the known solution vector $u_h^{k-1}\in V_h^{k-1}$ does not belong to the new solution space $V_h^k$. Before resolving the system on the new mesh, we simply recover (by interpolation or projection) the known solution vector $u_h^{k-1}$ to be used on the new solution space $V_h^k$.

\begin{algorithm}[tp]
\caption{Adaptive algorithm\label{adap_alg}}
Given initial mesh $\mathcal{T}_h^{0}$, initial space $V_h^{0}$ and initial vector $u_h^{0}\in V_h^{0}$:
\begin{algorithmic}[1]
\For{$k=1,2,\ldots , J$}
	\State  On the mesh $\mathcal{T}_h^{k-1}$, solve \eqref{dgAC} for $u_h^k\in V_h^{k-1}$
	\State $\forall E\in\mathcal{T}_h^{k-1}$, compute error indicators $\eta_E^k$
	\If{$(\eta_E^k)^2 < {\bf stol^r}$ \& $(\eta_E^k)^2 > {\bf stol^c}$, $\forall E$ }
		\State Set $\mathcal{T}_h^{k}:= \mathcal{T}_h^{k-1}$
		\State Set $V_h^{k}:= V_h^{k-1}$
	\Else
		\State Find $M_R, \; M_C\subset \mathcal{T}_h^{k-1}$ such that
		\State \quad $(\eta_E^k)^2 > {\bf stol^r} \Rightarrow E\in M_R$
		\State \quad $(\eta_E^k)^2 < {\bf stol^c} \Rightarrow E\in M_C$
		\State Form the new mesh $\mathcal{T}_h^{k}$, and new space $V_h^k$
		\State \quad Refine the elements $E\in M_R\subset\mathcal{T}_h^{k-1}$
		\State \quad Coarsen the elements $E\in M_C\subset\mathcal{T}_h^{k-1}$
		\State Interpolate (project) $u_h^{k-1}$ onto $V_h^{k}$
		\State On the mesh $\mathcal{T}_h^{k}$, resolve \eqref{dgAC} for $u_h^k\in V_h^k$
	\EndIf	
	\State Return $\mathcal{T}_h^{k}$, $V_h^{k}$ and $u_h^k$
\EndFor
\end{algorithmic}
\end{algorithm}

The accuracy of the solutions is through a better choice for the initial mesh $\mathcal{T}_{h}^{0}$ by pre-refining the uniform initial mesh so that the approximation error of the initial condition is well enough. Thus, in order to determine the initial mesh $\mathcal{T}_{h}^{0}$, we form a sequence of initial meshes $\{\mathcal{T}_{h}^{0,k}\}$ for $k=0,1,\ldots $, and we set the initial mesh $\mathcal{T}_{h}^{0}:=\mathcal{T}_{h}^{0,K}$ for some integer $K\geq 1$, where $\mathcal{T}_{h}^{0,0}$ is the initial uniform mesh and each $\mathcal{T}_{h}^{0,k}$ is obtained by refining $\mathcal{T}_{h}^{0,k-1}$. To select the elements to be refined in $\mathcal{T}_{h}^{0,k-1}$, we simply form the set
\[
M_{R,k} = \{ E\in \mathcal{T}_{h}^{0,k-1}\; : \; (\rho_E^{k-1})^2 > \bf{stol^0} \},
\]
where $\bf{stol^0}$ is a user defined tolerance, and $\rho_E^{k-1}$ is the approximation error given by
\[
(\rho_E^{k-1})^2 := \|g(\bm{x})-g_h(\bm{x})\|_{L^2(E)}, \quad E\in\mathcal{T}_{h}^{0,k-1},
\]
We select the initial mesh when it is satisfied that $\max_{E\in\mathcal{T}_{h}^{0,k}}\leq\bf{stol^0}$, or a maximum number $K$ of pre-refinement level is reached. The adaptive procedure is given in Algorithm~\ref{adap_alg}.

\section{Numerical results}
\label{sec:numerics}

In this section, we present numerical examples for advective AC equation under homogenous Neumann boundary conditions to demonstrate the effectiveness of the adaptive SIPG method to recapture sharp layers. All simulations are performed on the spatial domain $\Omega =[-1,1]^2$ using linear polynomials in the dG space. The final time is $T=0.06$, and the step-size is taken as $\tau = 0.001$. In the figures, we indicate the (spatial) dimension of dG space by DoFs which is  the product of the number of triangular elements in the mesh and the local dimension on the elements, i.e. the number of dG basis functions defined on each element.

\subsection{Expanding flow}
\label{ex1}

\begin{figure}[htp]
 \begin{center}
	\subfloat{\includegraphics[width=1.7in]{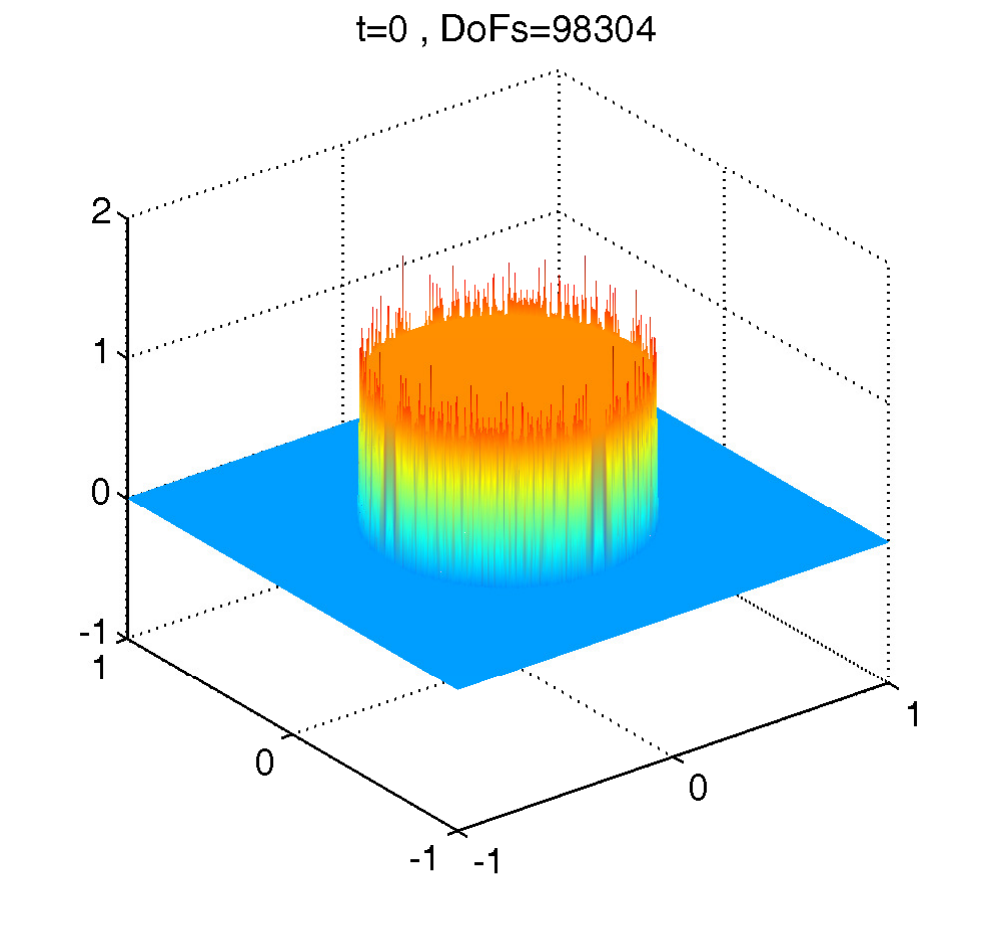}}
	\subfloat{\includegraphics[width=1.7in]{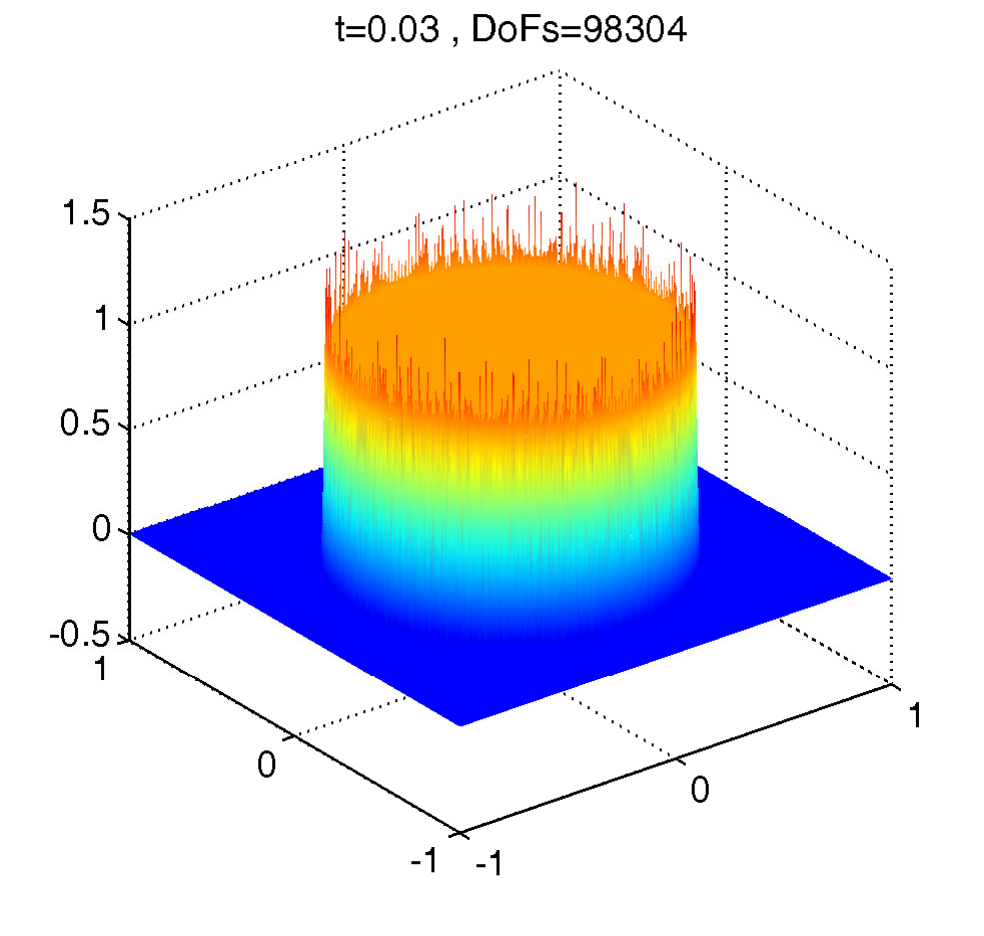}}
	\subfloat{\includegraphics[width=1.7in]{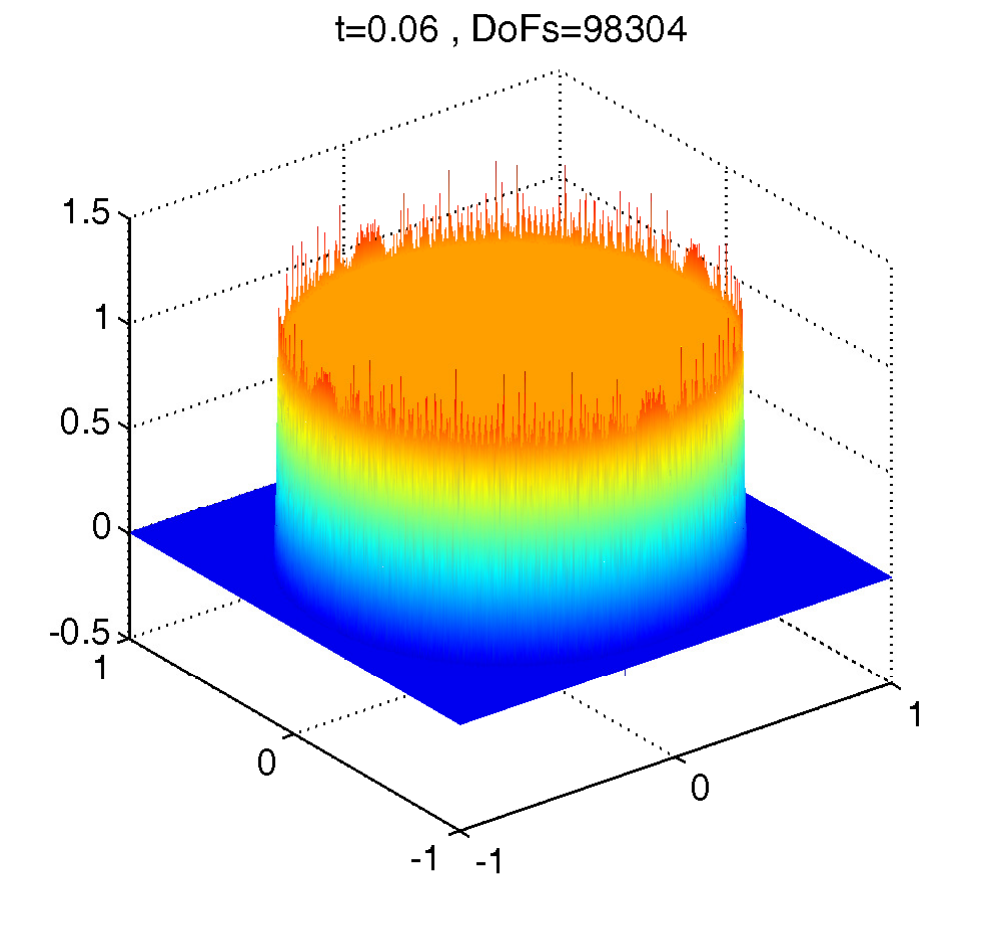}}
	
	\subfloat{\includegraphics[width=1.7in]{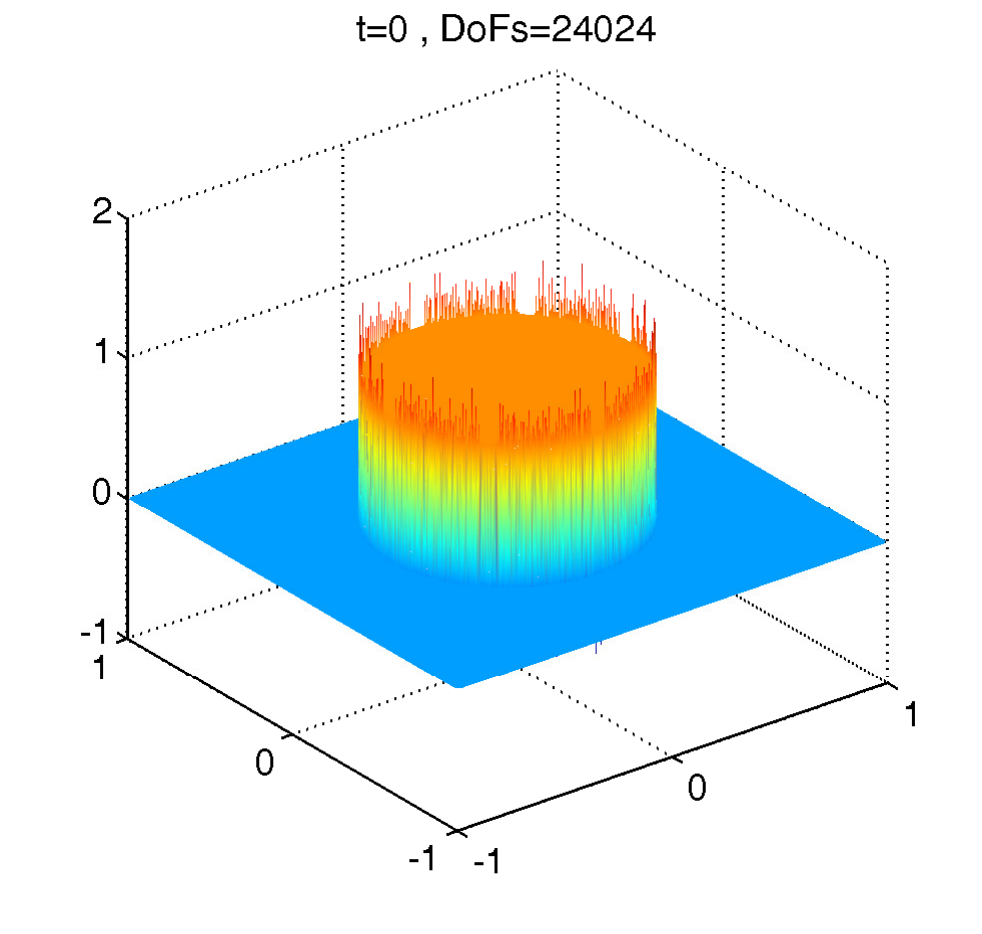}}
	\subfloat{\includegraphics[width=1.7in]{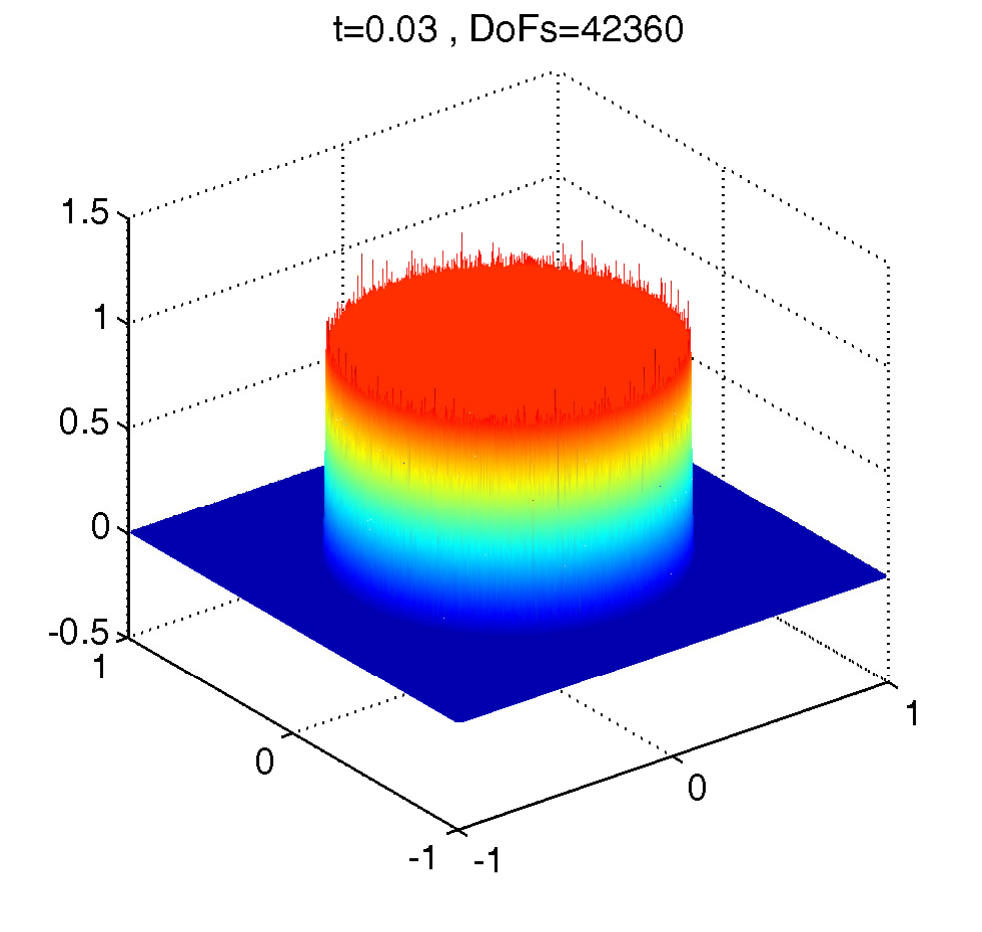}}
	\subfloat{\includegraphics[width=1.7in]{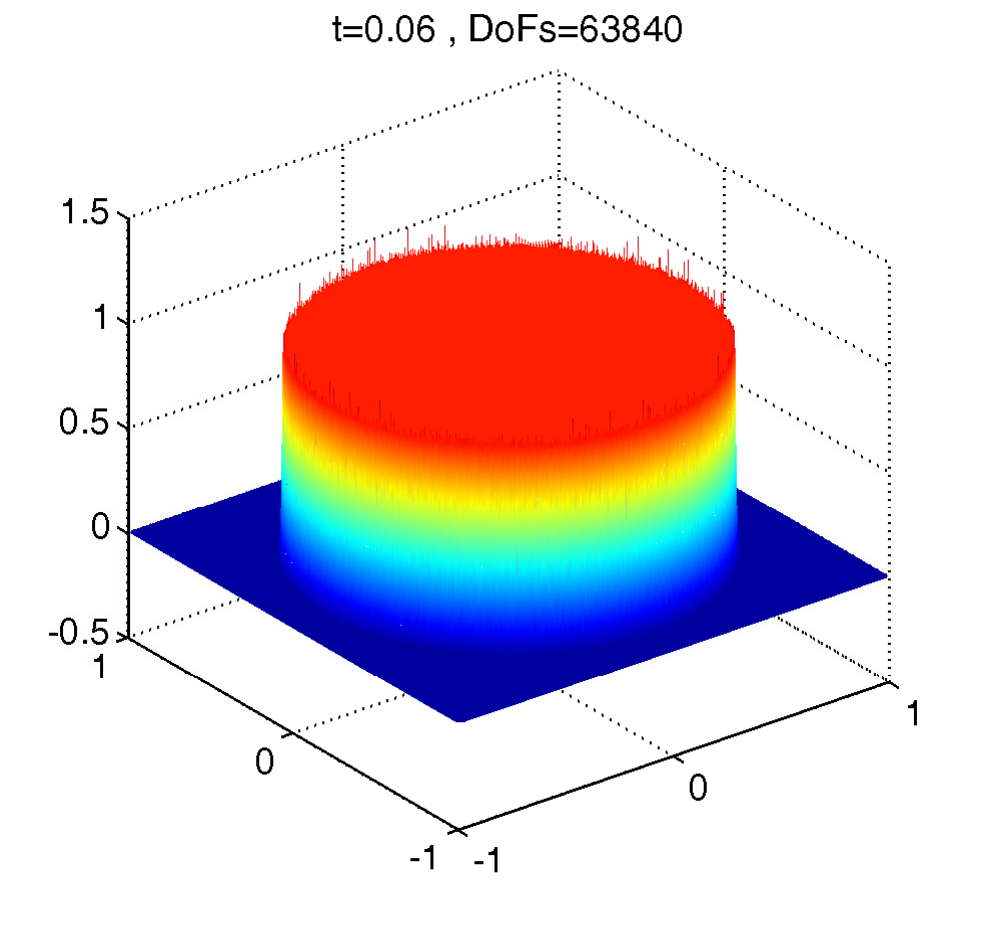}}
	
	\subfloat{\includegraphics[width=1.7in]{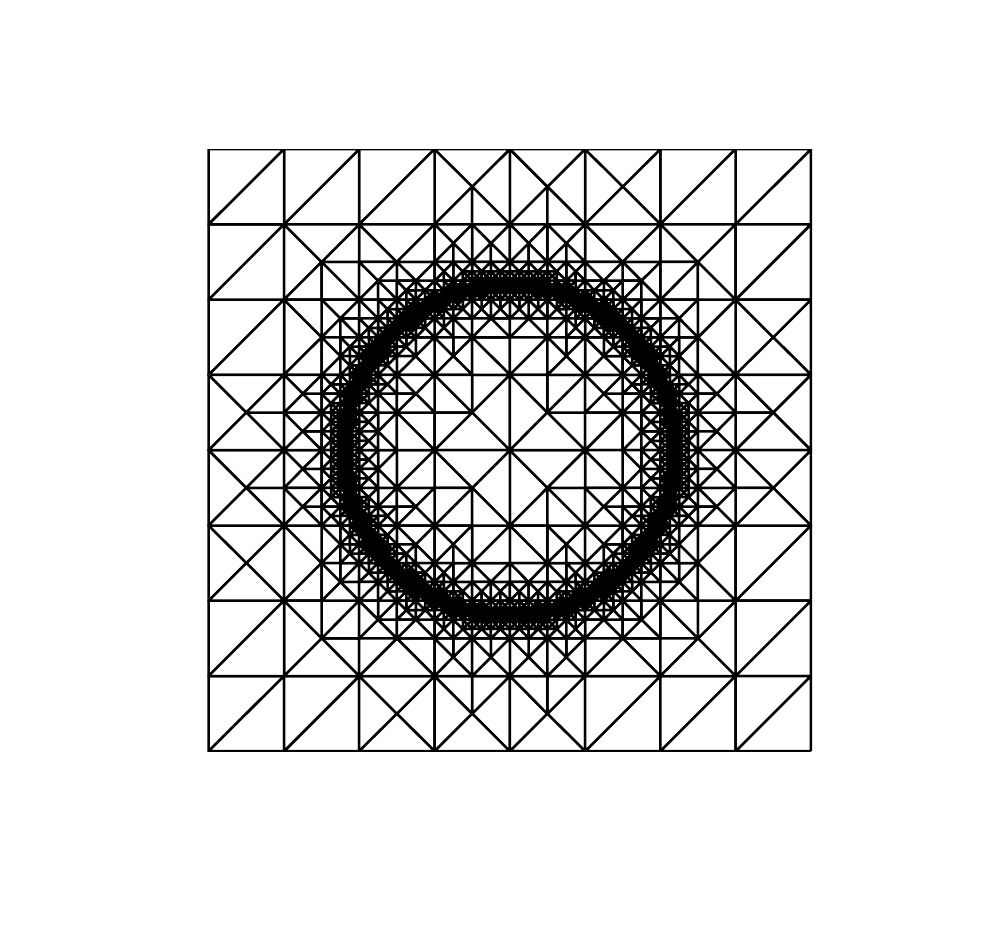}}
	\subfloat{\includegraphics[width=1.7in]{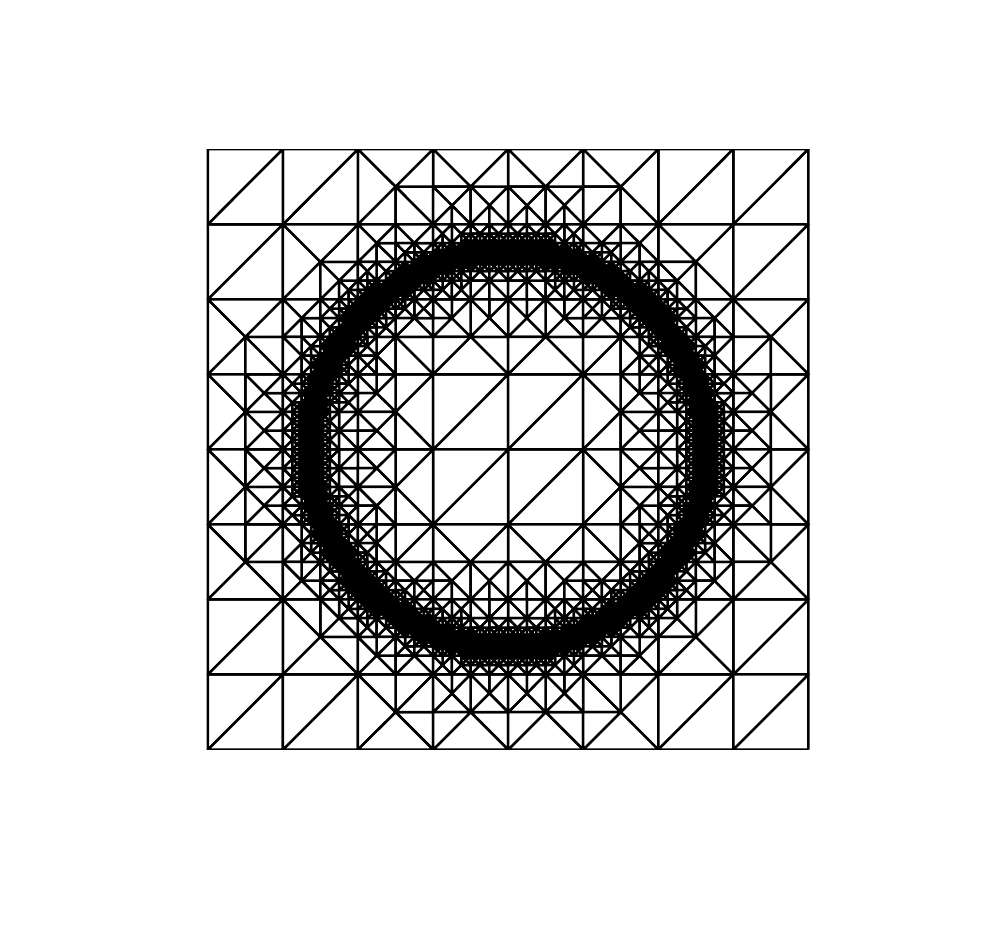}}
	\subfloat{\includegraphics[width=1.7in]{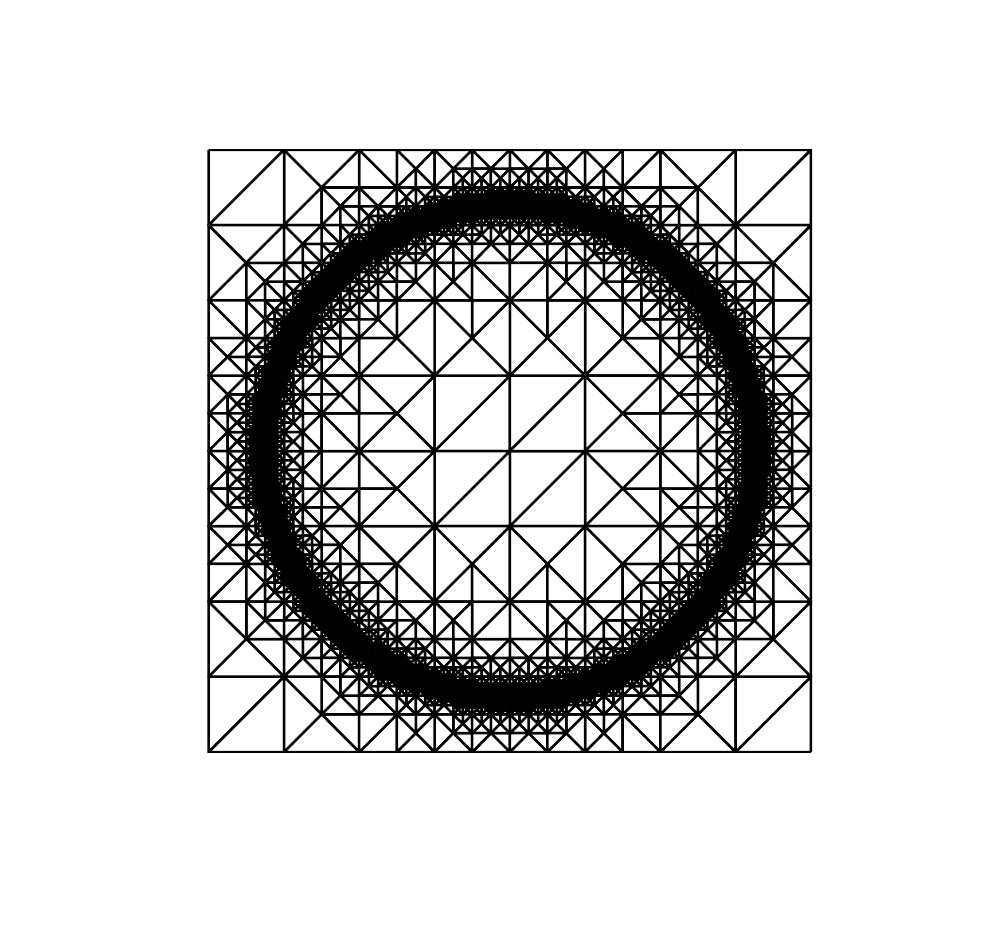}}
  \caption{Expanding flow: (Top) Uniform solutions, (middle) adaptive solutions and (bottom) adaptive meshes at time instances $t=0$, $t=0.03$ and $t=0.06$ from left to right}\label{sym_plot}
 \end{center}
\end{figure}

\begin{figure}[htp]
\begin{center}
	\includegraphics[width=2in] {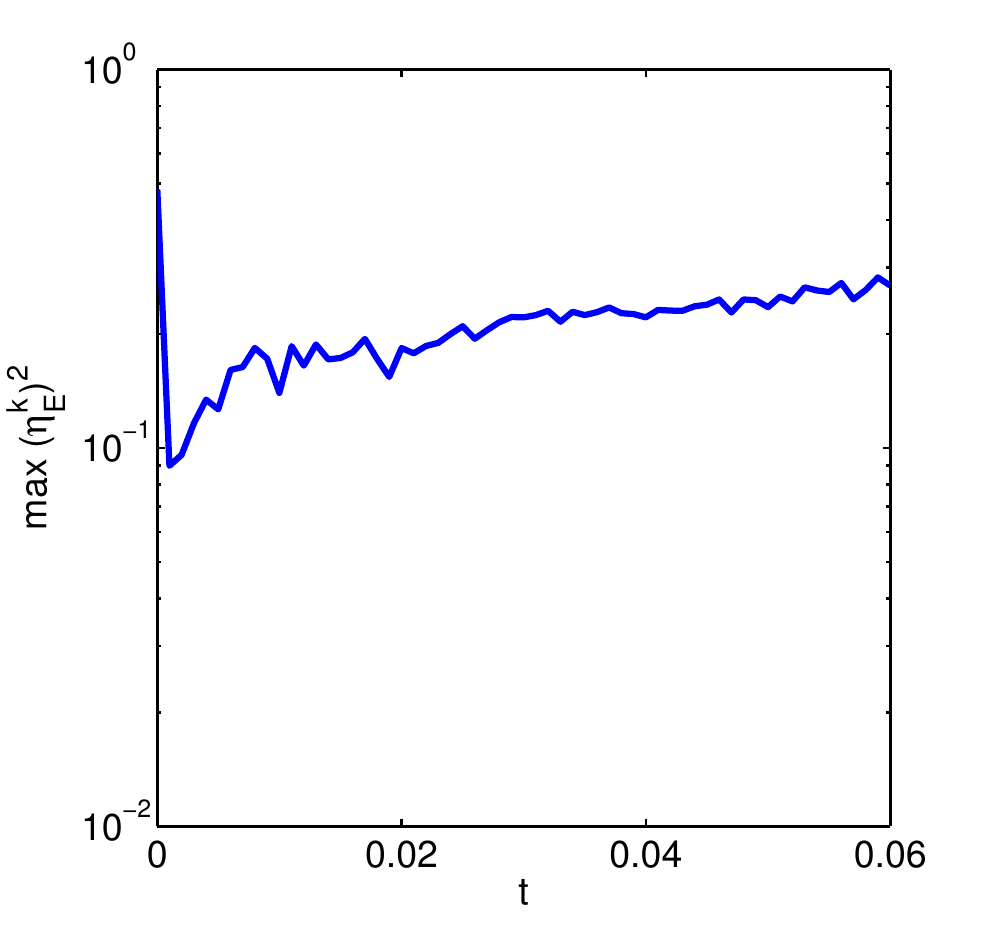}
	\includegraphics[width=2in] {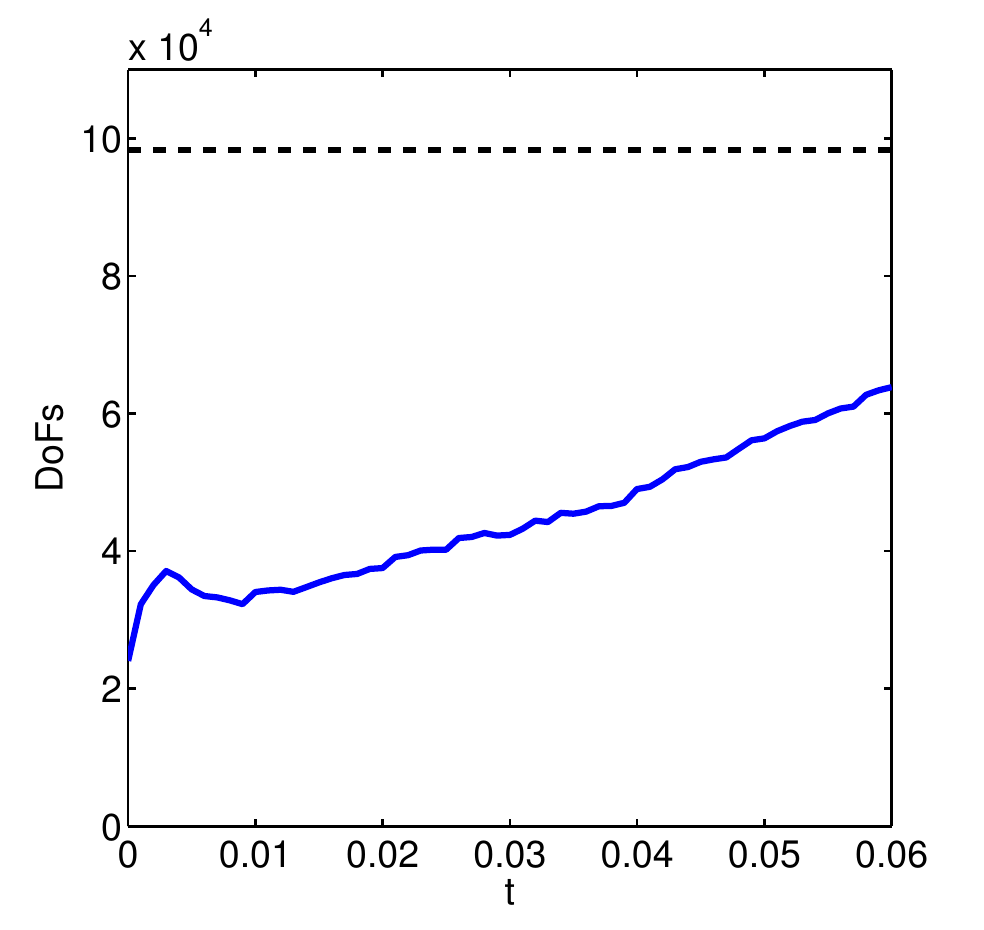}
\caption{Expanding flow: (Left) Maximum element error propagation over time, and (right) evaluation of DoFs over time; dashed line indicates the DoFs used for uniform solutions}\label{sym_error}
\end{center}
\end{figure}

We consider the advective AC equation with an expanding velocity field as \cite{liu12dis}
\[
\bm{V}=(v_0 x,v_0 y),
\]
with $v_0=10$. The initial condition is taken as the symmetric data
\[
u\left(x,0\right)=\begin{cases}
1;& x^2 + y^2 \leq 0.3\\
0;&\hbox{otherwise}
\end{cases},
\]
and the interface length is $\epsilon = 0.001$.
We first obtain the solutions on a uniform mesh with the mesh size $\Delta x=\Delta y=1/64$. Figure~\ref{sym_plot}, top, shows the uniform solutions at different time instances, where there seem spurious oscillations. For the adaptive case, we prescribe the tolerances $\bf{stol^r}=1\times 10^{-2}$ and  $\bf{stol^c}=1\times 10^{-5}$, and we take the initial uniform mesh with the mesh size $\Delta x=\Delta y=1/4$. The adaptive solutions and corresponding adaptive meshes are given in Figure~\ref{sym_plot}, middle-bottom. We can see that with much less DoFs, the oscillations are almost disappeared, and that the internal layers are well-captured together with the expanding behavior. In Figure~\ref{sym_error}, we present the propagation of the maximum element error $\eta_E^k$ over the time, and also the evaluation of DoFs during the time progression. We see that maximum element error is bounded in a small band according to the prescribed tolerances, and both refinement and coarsening phenomena works well.

\subsection{Sheer flow}
\label{ex2}

\begin{figure}[htp]
 \begin{center}
	\subfloat{\includegraphics[width=1.7in] {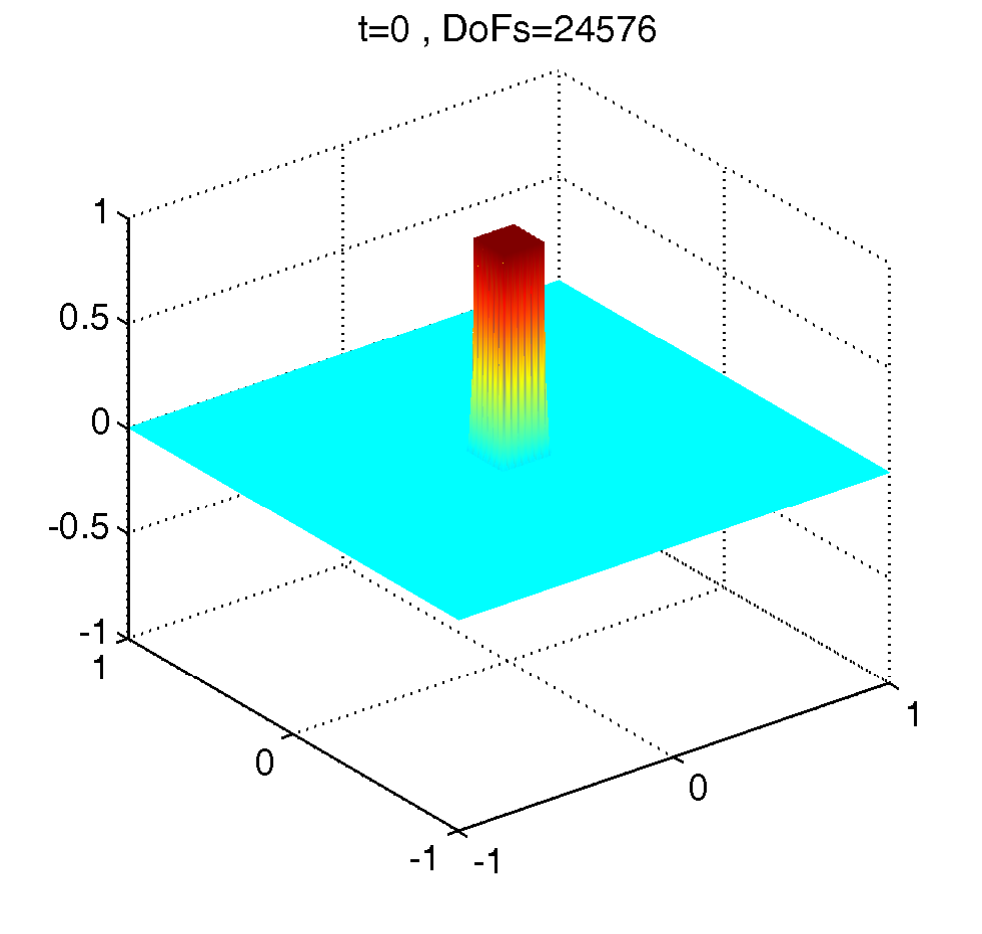}}
	\subfloat{\includegraphics[width=1.7in] {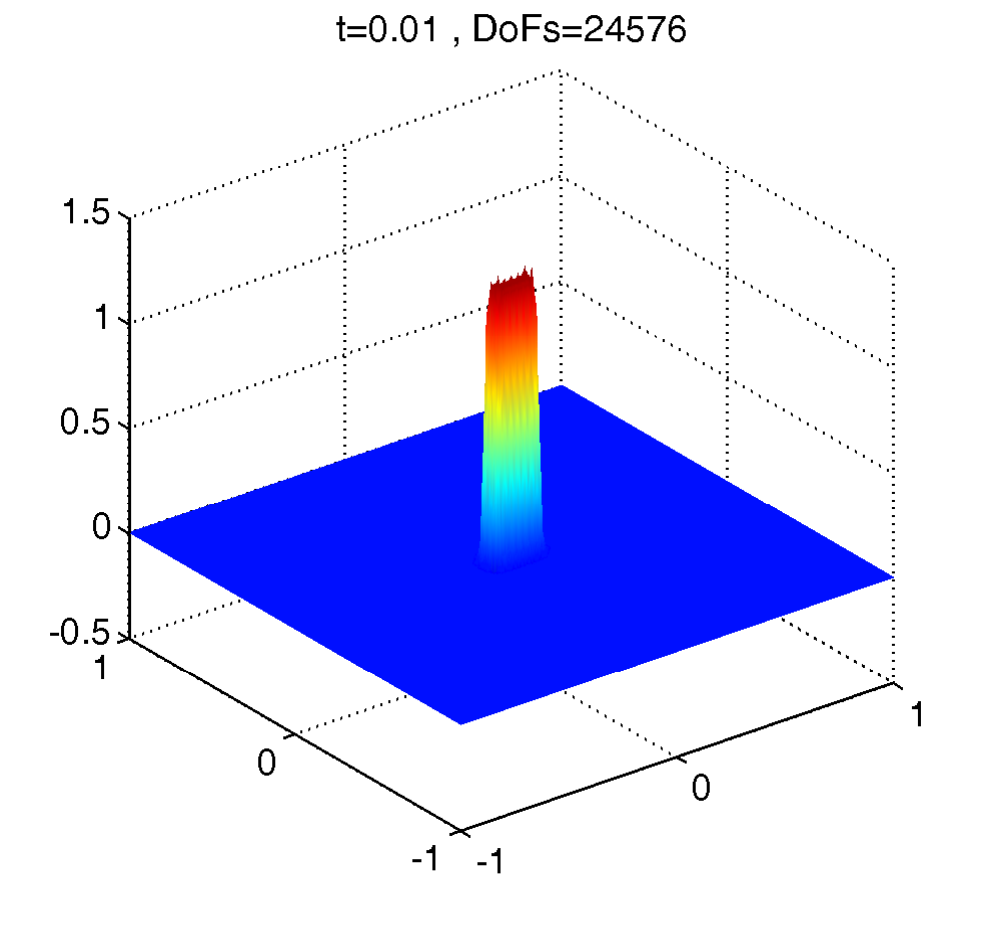}}
	\subfloat{\includegraphics[width=1.7in] {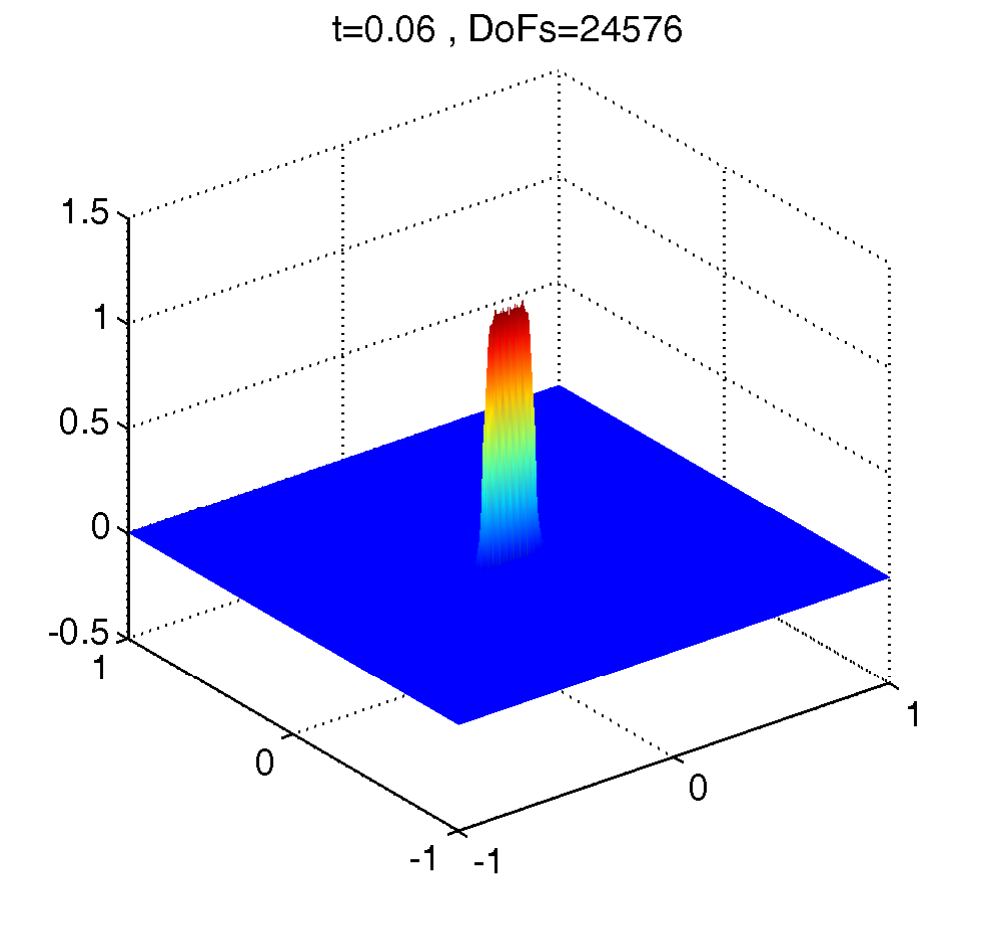}}
	
	\subfloat{\includegraphics[width=1.7in] {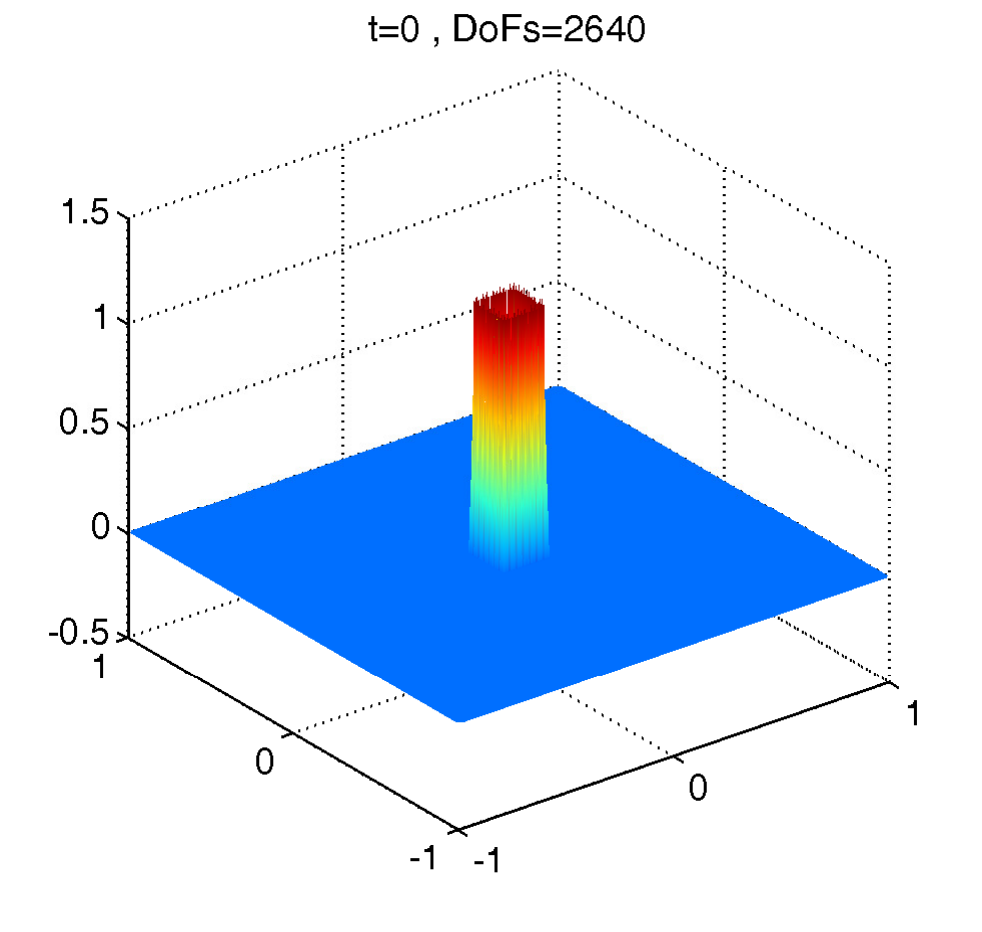}}
	\subfloat{\includegraphics[width=1.7in] {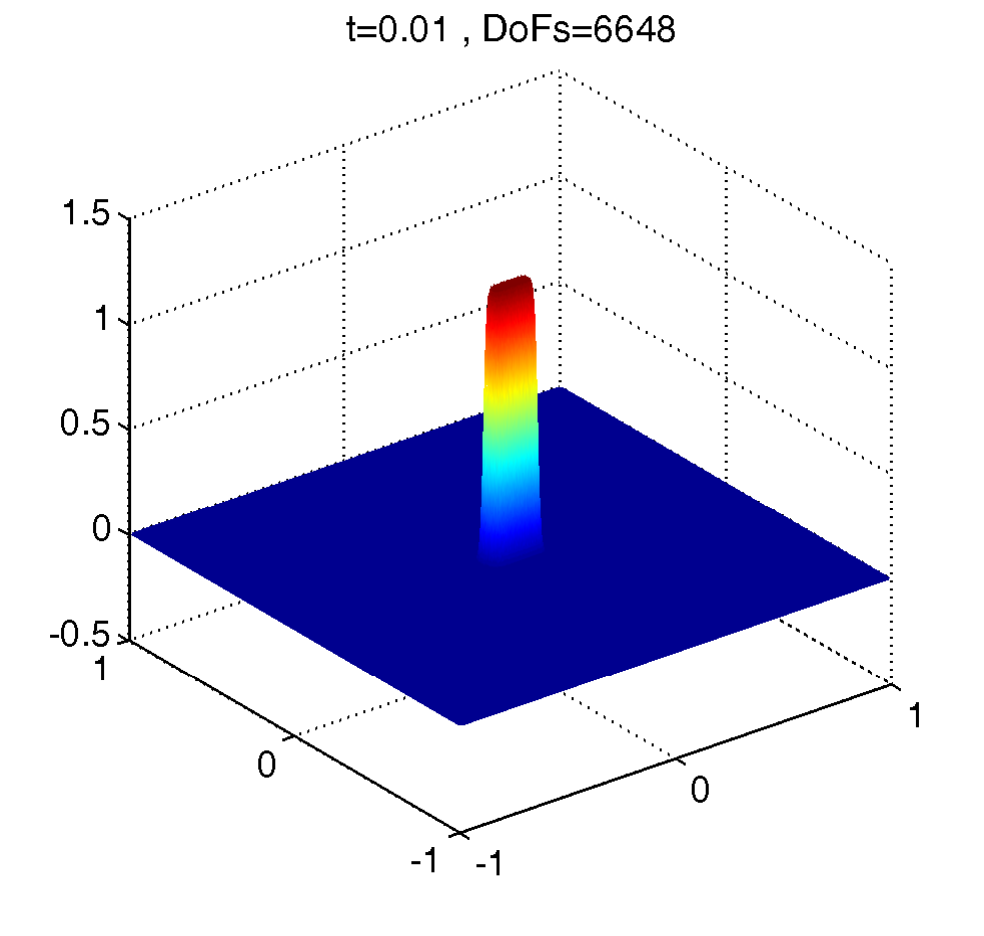}}
	\subfloat{\includegraphics[width=1.7in] {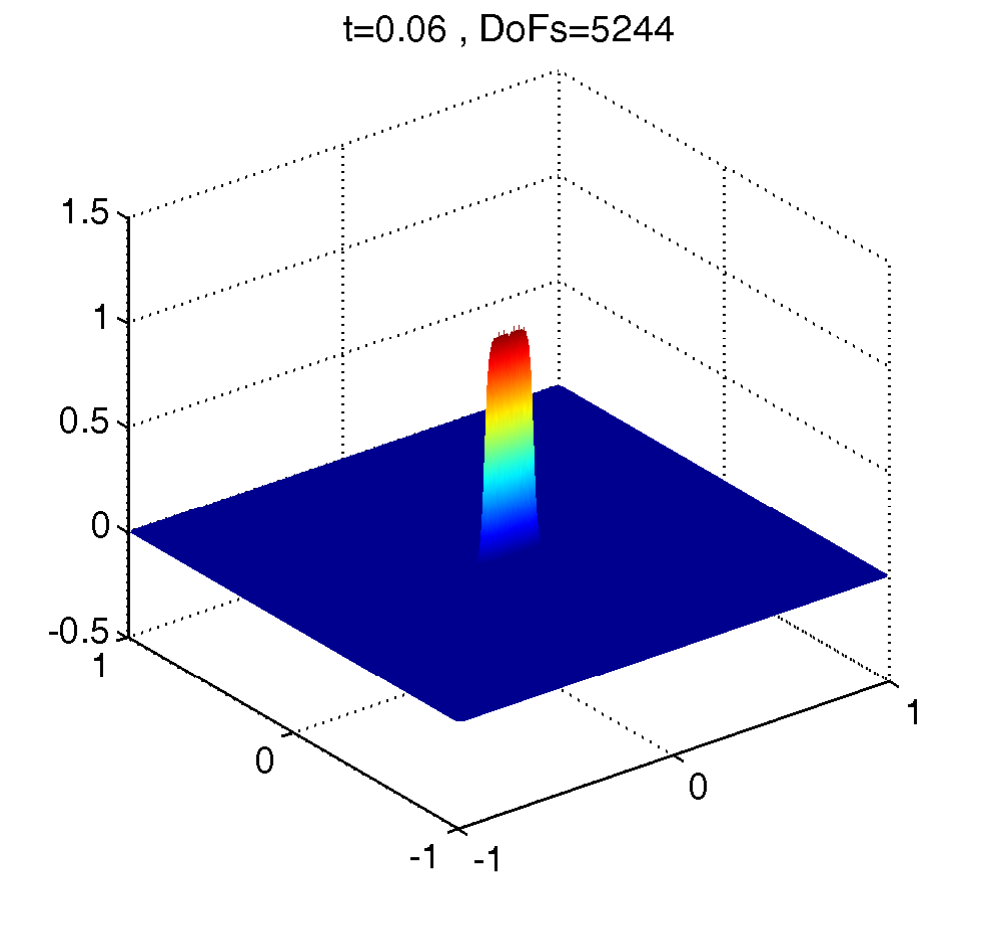}}
	
	\subfloat{\includegraphics[width=1.7in] {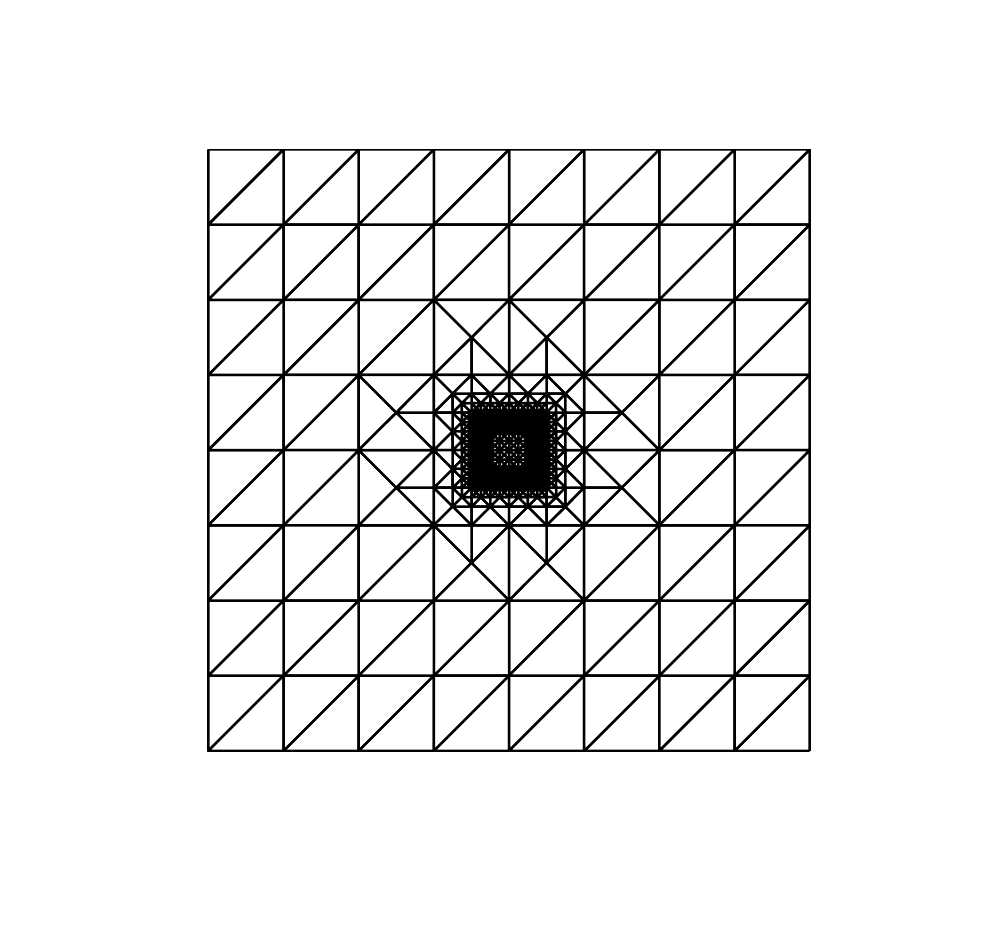}}
	\subfloat{\includegraphics[width=1.7in] {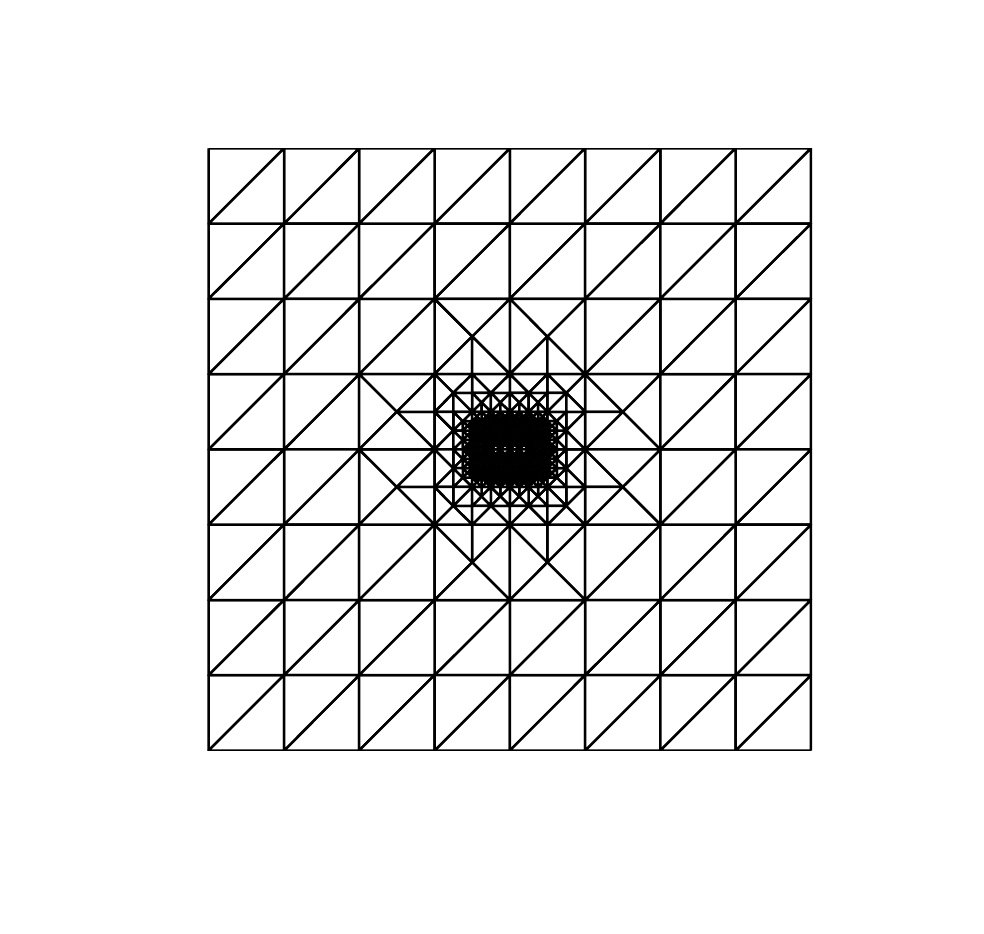}}
	\subfloat{\includegraphics[width=1.7in] {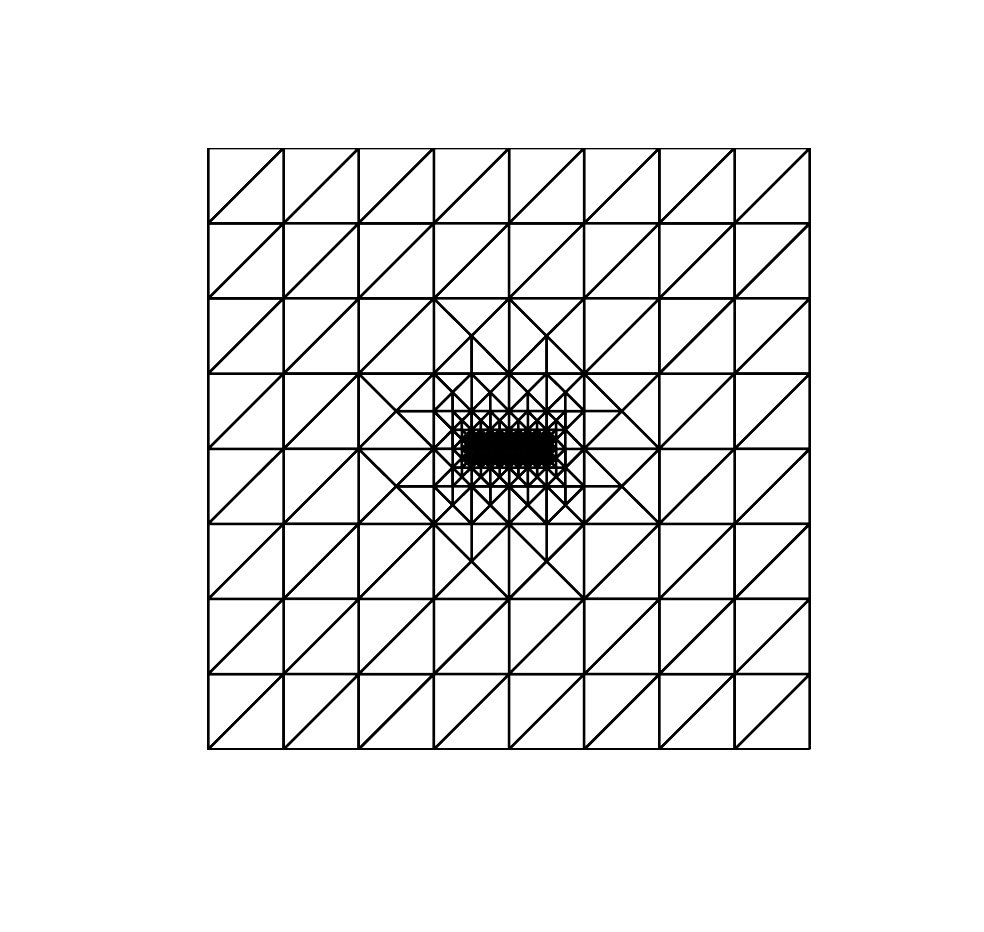}}
 \caption{Sheering flow: (Top) Uniform solutions, (middle) adaptive solutions and (bottom) adaptive meshes at time instances $t=0$, $t=0.01$ and $t=0.06$ from left to right }\label{sheer_plot}
 \end{center}
\end{figure}

\begin{figure}[htp]
 \begin{center}
	\includegraphics[width=2in]{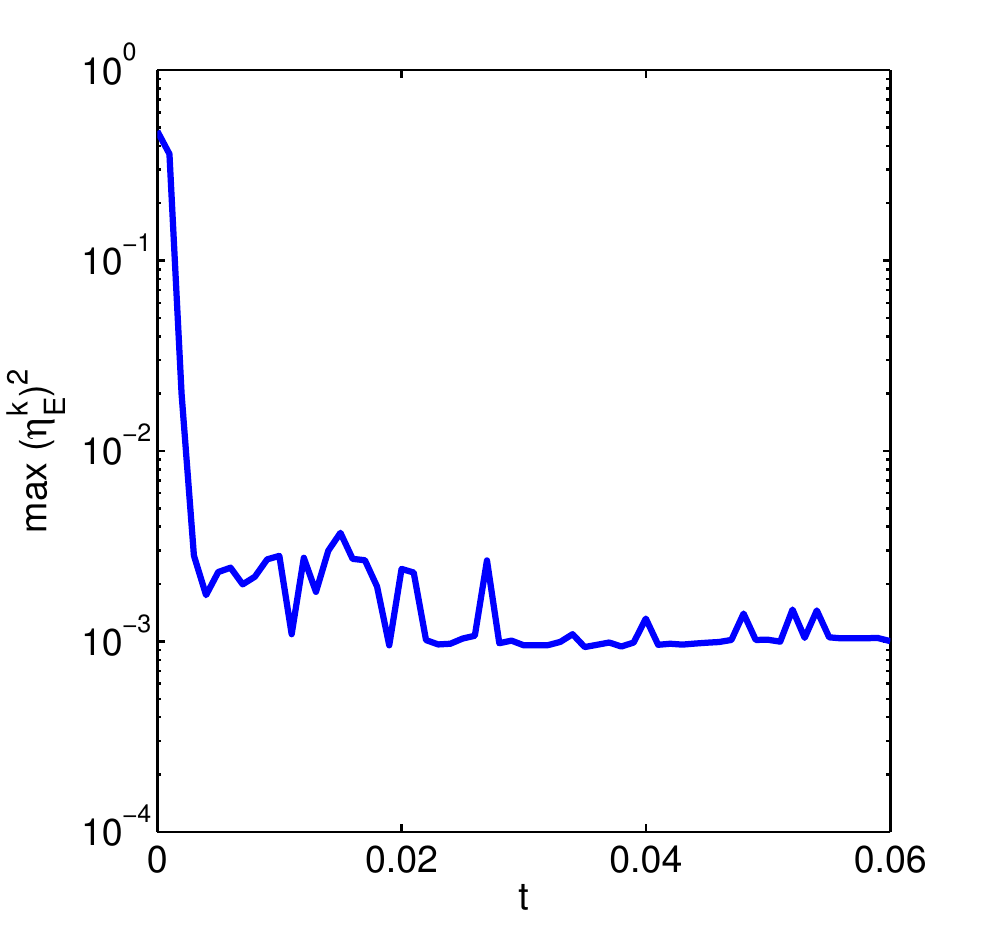}
	\includegraphics[width=2in]{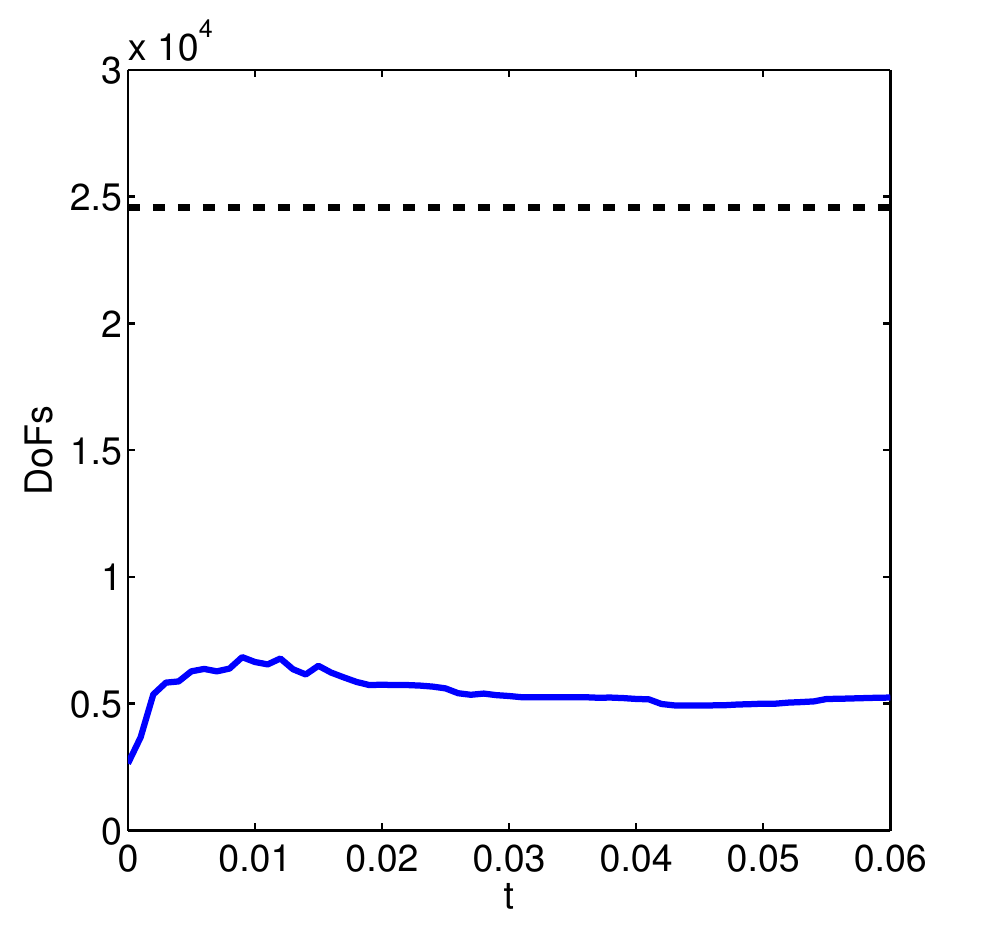}
  \caption{Sheering flow: (Left) Maximum element error propagation over time, and (right) evaluation of DoFs over time; dashed line indicates the DoFs used for uniform solutions}\label{sheer_error}
	\end{center}
\end{figure}

Now, we test the advective AC equation with a sheering flow given by \cite{liu12dis}
\[
\bm{V}=(0,-v_0 y),
\]
with $v_0=100$, and the initial condition is the square data
\[
u\left(x,0\right)=\begin{cases}
1;&-0.1\leq x, y\leq 0.1\\
0;&\hbox{otherwise}
\end{cases},
\]
The interface length is taken as $\epsilon = 0.01$. As the expanding case, the uniform solutions are obtained on the uniform mesh with the mesh size $\Delta x=\Delta y=1/32$, and for the adaptive case we use the same settings. In Figure~\ref{sheer_plot}, top, uniform solutions are shown with small oscillations, whereas they are dumped out in the adaptive solutions, Figure~\ref{sheer_plot}, middle. The layers are well-captured by the related adaptive meshes in Figure~\ref{sheer_plot}, bottom, in coherence with sheering behavior of the solutions. Moreover, the maximum element errors again lies in a small band, and the number of DoFs needed for the adaptive scheme is pretty less then it requires for the uniform case, Figure~\ref{sheer_error}.

\section{Conclusion}
Most of the adaptive methods are developed for linear and nonlinear diffusion-convection-reaction equations with non-divergence-free velocity fields. In this paper we have developed an adaptive procedure for AC equation in compressible fluids using interior penalty discontinuous Galerkin method.  Numerical results show that comparing with the uniform meshes, the sharp layers can be resolved accurately using less number of DoFs on adaptive meshes.

\section*{Acknowledgments} The authors would like to thank the reviewer for the comments and suggestions that helped to improve the manuscript.

\end{document}